\RequirePackage{pdf14}

\documentclass[final]{siamltex}

\usepackage[T1]{fontenc}
\UseRawInputEncoding

\usepackage{amssymb,latexsym,amsmath}
\usepackage{mathrsfs}
\usepackage{enumitem}
\usepackage{caption}
\usepackage{comment}
\usepackage{inputenc}
 \usepackage{makeidx}
\usepackage{color}

\allowdisplaybreaks

\newtheorem{assumption}{Assumption}[section]
\newtheorem{example}{Example}[section]
\newtheorem{remark}{Remark}[section]

\allowdisplaybreaks

\usepackage{epsfig} 
\usepackage{verbatim}

\hypersetup{colorlinks={true},linkcolor={blue},citecolor=green}

\usepackage{amssymb,nicefrac,bm,upgreek,mathtools,verbatim}
\usepackage{dsfont}
\usepackage{mathrsfs}

\RequirePackage[normalem]{ulem}

\usepackage{color}

\newcommand{\sym}[1]{{\color{black} #1}}
\allowdisplaybreaks

\def\qed{\hfill $\diamond$}

\newcommand{\ttup}[1]{\textup{(}#1\textup{)}}

\newcommand{\smid}{\,|\,}               
\newcommand{\df}{\coloneqq}             
\newcommand{\D}{\mathrm{d}}             
\newcommand{\RR}{\mathbb{R}}            
\newcommand{\NN}{\mathbb{N}}            

\newcommand{\Act}{\mathbb{U}}           
\newcommand{\Uadm}{\Gamma_{\mathsf{A}}} 
\newcommand{\Um}{\Gamma_{\mathsf{M}}}   
\newcommand{\Usm}{\Gamma_{\mathsf{S}}}  
\newcommand{\Usd}{\Gamma_{\mathsf{SD}}} 

\newcommand{\cG}{\mathcal{G}}           
\newcommand{\cH}{\mathcal{H}}           
\newcommand{\cT}{\mathcal{T}}           
\newcommand{\Bor}{{\mathcal{B}}}        
\newcommand{\Pm}{{\mathcal{P}}}         
\newcommand{\cM}{{\mathcal{M}}}         

\newcommand{\XX}{{\mathbb{X}}}          
\newcommand{\KK}{\mathbb{K}}            
\newcommand{\HH}{\mathbb{H}}            
\newcommand{\cU}{{\mathcal{U}}}         

\definecolor{dmagenta}{rgb}{.4,.1,.5}
\definecolor{dblue}{rgb}{.0,.0,.5}
\definecolor{mblue}{rgb}{.0,.0,.7}
\definecolor{ddblue}{rgb}{.0,.0,.4}
\definecolor{dred}{rgb}{.7,.0,.0}
\definecolor{dgreen}{rgb}{.0,.5,.0}
\definecolor{Eeom}{rgb}{.0,.0,.5}

\allowdisplaybreaks

\begin{document}

\title{On Borkar and Young Relaxed Control Topologies and Continuous Dependence of Invariant Measures on Control Policy \thanks{This research was partially supported by the Natural Sciences and Engineering Research Council of Canada (NSERC).}}



\author{Serdar Y\"{u}ksel \thanks{S. Y\"uksel is with the Dept. of Mathematics and Statistics, Queen's University, Kingston K7L 3N6, ON, Canada, (\email{yuksel@queensu.ca}) }}

%



\maketitle

\begin{abstract}
In deterministic and stochastic control theory, relaxed or randomized control policies allow for versatile mathematical analysis (on continuity, compactness, convexity and approximations) to be applicable with no artificial restrictions on the classes of control policies considered, leading to very general existence results on optimal measurable policies under various setups and information structures. On relaxed controls, two studied topologies are the Young and Borkar (weak$^*$) topologies on spaces of functions from a state/measurement space to the space of probability measures on control action spaces; the former via a weak convergence topology on probability measures on a product space with a fixed marginal on the input (state) space, and the latter via a weak$^*$ topology on randomized policies viewed as maps from states/measurements to the space of signed measures with bounded variation. We establish implication and equivalence conditions between the Young and Borkar topologies on control policies. We then show that, under some conditions, for a controlled Markov chain with standard Borel spaces the invariant measure is weakly continuous on the space of stationary control policies defined by either of these topologies. An implication is near optimality of quantized stationary policies in state and actions or continuous stationary and deterministic policies for average cost control under two sets of continuity conditions (with either weak continuity in the state-action pair or strong continuity in the action for each state) on transition kernels.


\end{abstract}



\section{Introduction}


In deterministic and stochastic control theory, relaxed or randomized control policies allow for versatility in mathematical analysis leading to continuity, compactness, convexity and approximation properties, in a variety of system models, cost criteria, and information structures. 



Under the relaxed/randomized control framework, with $\mathbb{X}$ a standard Borel state space, $\mathbb{U}$ a standard Borel control space and with an $\mathbb{X}$-valued random variable $X \sim \mu$, instead of considering the set of deterministic admissible policies:
\begin{align}
\Gamma = \bigg\{\gamma: \text{$\gamma$ is a measurable function from $\mathbb{X}$ to $\mathbb{U}$} \bigg\}, \label{detPolS}
\end{align}
one considers
\begin{align}
\Gamma_R = \bigg\{\gamma: \text{$\gamma$ is a measurable function from $\mathbb{X}$ to ${\cal P}(\mathbb{U})$} \bigg\}, \label{ranPolS}
\end{align}
where ${\cal P}(\mathbb{U})$ is endowed with the Borel $\sigma$-algebra generated by the weak convergence topology. 

On $\Gamma_R$, two commonly studied topologies are the following.\\

\noindent{\bf Young Topology on randomized policies.} A particularly prominent approach since Young's seminal paper \cite{young1937generalized} has been via the study of topologies on {\it Young measures} defined by randomized/relaxed controls, where one views policies to be identified with probability measures defined on a product space with a fixed marginal at an input/state space (typically taken to be the Lebesgue measure in optimal deterministic control) \cite{young1937generalized,mcshane1967relaxed},\cite[Section 2.1]{castaing2004young},\cite[p. 254]{warga2014optimal},\cite{mascolo1989relaxation}, \cite[Theorem 2.2]{balder1988generalized}. 

Thus, under the Young topology, one associates with $\Gamma_R$ in (\ref{ranPolS}) the probability measure induced on the product space $\mathbb{X} \times \mathbb{U}$ with a fixed marginal $\mu$ on $\mathbb{X}$. \sym{On this product space, several weak topologies and their equivalence properties have been studied in the literature; e.g, \cite{balder1997consequences}, \cite[Theorem 2.2]{balder1988generalized} and \cite[Theorem 3.32]{florescu2012young} study the equivalence relations between the Young topology convergence and the narrow topology convergence}.

The generalization to stochastic control problems by considering more general input measures has been commonplace, with applications also to partially observed stochastic control and decentralized stochastic control as noted above, where most of the aforementioned references have indeed adopted this approach to define control topologies (see e.g. on existence and approximation results involving optimal stochastic control \cite{fleming1976generalized,kushner2001numerical,kushner2014partial}, on decentralized stochastic control \cite[Theorem 5.2]{YukselWitsenStandardArXiv} and \cite[Section 4]{saldiyukselGeoInfoStructure}, piece-wise deterministic optimal stochastic control \cite{bauerle2010optimal,bauerle2018optimal}, on economics and game theory \cite{milgrom1985distributional,mertens2015repeated,balder1988generalized}, mean-field control policies \cite{bayraktar2022finite} and on optimal quantization \cite[Definition 2.1]{YukselOptimizationofChannels}).\\

 \noindent{\bf Borkar's (weak$^*$) topology on randomized policies.} In the stochastic setup, another topology is one introduced by Borkar on relaxed controls \cite{borkar1989topology} (see also \cite[Section 2.4]{arapostathis2012ergodic}, and \cite{bismut1973theorie} which \cite{borkar1989topology} notes to be building on), formulated as a weak$^*$ topology on randomized policies viewed as maps from states/measurements to the space of signed measures with bounded variation ${\cal M}(\mathbb{U})$ of which probability measures ${\cal P}(\mathbb{U})$ is a subset. We also refer the reader to \cite{dieudonne1951theoreme,fattorini1994existence} for further references on such a weak$^*$ formulation on relaxed controls, in particular when instead of countably additive signed measures, finitely additive such measures are also considered.

See \cite{borkar1989topology,arapostathis2010uniform,pradhan2022near} for a detailed analysis on some implications in stochastic control theory in continuous-time (such as continuity of expected cost in control policies \cite{borkar1989topology}, approximation results \cite{pradhan2022near} under various cost criteria, and continuity of invariant measures of diffusions in control policies \cite{arapostathis2010uniform}).  

Thus, under the Borkar topology one studies $\Gamma_R$ in (\ref{ranPolS}), with a weak$^*$ topology formulation, as a bounded subset of the set of maps from $\mathbb{X}$ to the space of signed measures with finite variation viewed as the topological dual of continuous functions vanishing at infinity, leading to a compact metric space by the Banach-Alaoglu theorem \cite[Theorem 5.18]{Fol99} (and thus, as the unit ball of $L_{\infty}(\mathbb{X}, {\cal M}(\mathbb{U})) = (L_{1}(\mathbb{X}, C_0(\mathbb{U})))^*$ is compact under the weak$^*$ topology, this leads to a compact metric topology on relaxed control policies). We note that the presentations in \cite[Section 3]{borkar1989topology} and \cite[Section 2.4]{arapostathis2012ergodic} are slightly different, though the induced topology is identical. An equivalent representation of this topology is given in \cite[Lemma 2.4.1]{arapostathis2012ergodic} (see also \cite[Lemma 3.1]{borkar1989topology}). 

While $\mathbb{X}$ was taken to be $\mathbb{R}^n$ in \cite{borkar1989topology}, Saldi generalized this to setups where $\mathbb{X}$ is a general standard Borel space with a fixed input (probability) measure. This topology was used for optimal decentralized stochastic control for an existence analysis by Saldi \cite{SaldiArXiv2017}.

Further discussion and relations on these topologies will be presented in the paper. 
 
It should also be noted that, while some further classical contributions on stochastic control theory, such as \cite{benevs1971existence,duncan1971solutions,bismut1982partially}, did not consider randomized policies for the analysis on optimal stochastic control problems, they did consider topologies which are closely related to the Young topology.\\

We will also consider an application of the above on continuity properties of invariant measures on the space of stationary control policies, where we will obtain a discrete-time counterpart of \cite[Lemma 4.4]{arapostathis2010uniform} (under the Young topology and for more general state spaces). This result, in addition to applications on existence and approximations of optimal control under average cost criterion which we will study in this paper, also finds direct applications in stochastic game theory (see e.g. \cite[Assumption A6]{biswas2015mean} for which our results will present sufficient conditions for a study of mean-field games with ergodic cost for a collection of agents coupled only via their cost and with local information, facilitating an existence analysis through the Kakutani-Fan-Glicksberg fixed point theorem). In addition to these applications, another motivation for comparing the topologies has been the following: Young topology was utilized in \cite{YukselWitsenStandardArXiv} to arrive at general existence and approximation results in decentralized stochastic control, whereas the the generalization of Borkar topology by Saldi was utilized to study the same in \cite{SaldiArXiv2017}. The apparent connection, and the continuity problem noted above, motivate the analysis in the paper.



{\bf Contributions.} Our main results are the following.
\begin{itemize}
\item[(i)] Building on a supporting technical result on the relations between Young topologies at absolutely continuous input measures in Lemma \ref{equivYoungT}, we provide in Theorem \ref{equivYouBorTop} implication and equivalence relations between the Borkar and Young topologies. Notably, with $\psi$ being an input measure and $\lambda$ the Lebesgue measure and $\mathbb{X} = \mathbb{R}^n$; if $\psi \ll \lambda$ with the Radon-Nikodym derivative (density) $g(x)=\frac{d\psi}{d\lambda}(x)$ being positive everywhere, convergence in Young topology at input measure $\psi$ implies convergence in Borkar topology; and if $\psi \ll \lambda$ and $\psi$ is a finite measure, then convergence in Borkar topology implies convergence in Young topology at input $\psi$. Thus, if $\psi \ll \lambda$, $\psi$ is a finite measure (such as a probability measure), and $g(x)=\frac{d\psi}{d\lambda}(x)$ is positive everywhere, the topologies are equivalent. 


\item[(ii)] We present conditions on the continuous dependence of invariant measures in discrete-time controlled Markov chains on the space of control policies under the Young topology in Theorem \ref{ContInvPolicy} and a particular application for finite models in Theorem \ref{ContInvPolicy2}. Theorem \ref{ContInvPolicy} has direct implications on existence of average cost optimal stochastic control policies (we recall that in continuous-time, such a result has been established \cite[Lemma 4.4]{arapostathis2010uniform} under the Borkar topology). In discrete-time, such a result is not available to our knowledge for controlled models. In a control-free setup, we note a perturbation result given in \cite[Thm. 1.6]{kartashov1996strong} (see \cite[Theorem 9.2.6]{HernandezLasserreErgodic}) under a strong norm defined on transition kernels viewed as an operator acting on the space of signed measures.
\item[(iii)] As a further application, in Theorem \ref{denseQuantizedFiniteAC}, we show near optimality of quantized control policies (in both state and action) or continuous stationary and deterministic policies for average cost control under mild technical (with either weak continuity in the state-action pair or strong continuity in the action for each state) conditions on transition kernels. These justify the use of learning theoretic methods or numerical methods for average cost optimal control problems. 
\end{itemize}

%
%

\section{Model}

We consider the usual model in the literature for controlled
Markov chains, otherwise referred to as Markov decision processes (MDPs).
In general, for a topological space $\mathcal{X}$, we denote by $\Bor(\mathcal{X})$
its Borel $\sigma$-field and by $\Pm(\mathcal{X})$ the set of probability
measures on $\Bor(\mathcal{X})$.

A \emph{controlled Markov chain} consists of the tuple $\bigl(\XX,\Act,\cU,\cT,c\bigr)$,
whose elements can be described as follows.
\begin{itemize}
\item[(a)]
The \emph{state space} $\XX$ and the
\emph{action} or
\emph{control space} $\Act$ are Borel subsets of complete, separable,
metric (i.e., \emph{Polish}) spaces.
\item[(b)]
The map $\cU\colon \XX\to \Bor(\Act)$ is a strict, measurable
multifunction.
The set of admissible state/action pairs is 
\begin{equation*}
\KK \,\df\, \bigl\{(x,u)\colon\, x\in\XX,\,u\in\cU(x)\bigr\}\,,
\end{equation*}
endowed with the subspace topology
corresponding to $\Bor(\XX\times\Act)$.
\item[(c)]
The map $\cT\colon\KK\to\Pm(\XX)$ is a stochastic kernel on
$\KK\times\Bor(\XX)$, that is, $\cT(\,\cdot\smid  x,u)$ is
a probability measure on $\Bor(\XX)$ for each $(x,u)\in\KK$,
and $(x,u) \mapsto \cT(A\smid  x,u)$ is measurable for each $A\in\Bor(\XX)$.
\end{itemize}

The (admissible) \emph{history spaces} are defined as
\begin{equation*}
\HH_0\,\df\, \XX\,, \quad \HH_{t} \,\df\, \HH_{t-1}\times \mathbb{U} \times \mathbb{X} ,\quad t\in\NN\,,
\end{equation*}
and the canonical sample space is defined as $\Omega\df (\XX\times\Act)^\infty$.
These spaces are endowed with their respective product topologies and
are therefore Borel spaces.
The state, action (or control), and information processes, denoted
by $\{X_t\}_{t\in\NN_0}$,  $\{U_t\}_{t\in\NN_0}$ and $\{H_t\}_{t\in\NN_0}$,
respectively, are defined by the projections
\begin{equation*}
X_t(\omega) \,\df\, x_t\,,\quad U_t(\omega) \,\df\, u_t\,,
\quad H_t(\omega) \,\df\, (x_0,u_0, \dotsc, u_{t-1}, x_t)
\end{equation*}
for each $\omega=(x_0,u_0, \dotsc,u_{t-1},x_t, u_t,\dotsc)\in\Omega$.
An \emph{admissible control policy}, or \emph{policy}, is a sequence
$\gamma = \{\gamma_t\}_{t\in\NN_0}$ of stochastic kernels on
$\HH_t\times\Bor(\Act)$ satisfying the constraint
\begin{equation*}
\gamma_t(\cU(X_t)\mid h_t) \,=\, 1\,,\quad  h_t\in\HH_t\,.
\end{equation*}
The set of all admissible policies is denoted by $\Uadm$.
It is well known
(see \cite[Prop.\ V.1.1, pp.~162--164]{Neveu})
that for any given $\nu\in\Pm(\XX)$ and $\gamma\in\Uadm$ there exists
a unique probability measure $P^\gamma_\nu$ on $\bigl(\Omega,\Bor(\Omega)\bigr)$
satisfying
\begin{align*}
P^\gamma_\nu(X_0\in D) &\,=\, \nu(D)\qquad\forall\, D\in\Bor(\XX)\,,\\
P^\gamma_\nu(U_t\in C\mid H_t) &\,=\, 
\gamma_t(C\mid H_t)
\quad P^\gamma_\nu\text{-a.s.}\,,\quad \forall\, C\in\Bor(\Act)\\
P^\gamma_\nu(X_{t+1}\in D\mid H_t,U_t)  &\,=\,
\cT(D\mid X_t, U_t) \quad P^\gamma_\nu\text{-a.s.}\,,
\quad\forall\, D\in\Bor(\XX)\,.
\end{align*}
The expectation operator corresponding to $P^\gamma_\nu$
is denoted by $E^\gamma_\nu$.
If $\nu$ is a Dirac mass at $x\in\XX$, we simply write these
as $P^\gamma_x$ and $E^\gamma_x$.

A policy $\gamma$ is called \emph{Markov}
if there exists a sequence of measurable maps $\{v_t\}_{t\in\NN_0}$,
where $v_t\colon \XX\to \Pm(\Act)$, where $\Pm(\Act)$ is endowed with the weak convergence topology, for each $t\in\NN_0$,
such that
$$\gamma_t(\,\cdot\mid H_t) \,=\, v_t(X_t)(\cdot) \quad P^\gamma_\nu\text{-a.s.}$$
With some abuse in notation, such a policy is identified
with the sequence $v=\{v_t\}_{t\in\NN_0}$.
Note then that
 $\gamma_t$ may be written as a stochastic kernel $\gamma_t(\cdot\smid  x)$ on
$\XX\times\Bor(\Act)$ which satisfies $\gamma_t(\cU(x)\smid  x)=1$.
Let $\Um$ denote the set of all Markov policies. 


We add the adjective \emph{stationary} to indicate that the Markov policy
does not depend on $t\in\NN_0$, that is, $\gamma_t=\gamma$ for all $t\in\NN_0$.
We let $\Usm$ denote the class of stationary Markov policies,
henceforth referred to simply as \emph{stationary policies},
and let $\Usd\subset\Usm$ denote the subset of those that are deterministic.

In summary, under a policy $\gamma\in\Usm$, the process
satisfies the following: for all Borel sets $B \in \mathcal{B}(\XX)$, $t \ge 0$, and
($P^\gamma$ almost all) realizations $x_{[0,t]}, u_{[0,t]}$, we have
\begin{equation} \label{eq_evol}
\begin{aligned}
P^\gamma\bigl( X_{t+1} \in B \smid  X_{[0,t]}=x_{[0,t]}, U_{[0,t]}=u_{[0,t]}\bigr)
&\,=\, P^\gamma( X_{t+1} \in B \smid  X_t=x_t, U_t=u_t) \\
&\,= \cT(  B \smid  x_t, u_t)\,.
\end{aligned}
\end{equation}



For $\gamma\in\Usm$, we let
\begin{equation}\label{E-Tgamma}
\cT^\gamma (A\smid x) \,\df\, \int_{\cU(x)}
\cT(A\smid x,u)\,\gamma(\D{u}\smid  x)\,.
\end{equation}

We let $\cM_b(\XX)$ ($C_b(\XX)$)
 denote the space of bounded Borel measurable (continuous) real-valued functions
on $\XX$.
For $\mu\in\Pm(\KK)$ and $f\in\cM_b(\XX)$, we define
$\mu\cT\in\Pm(\XX)$ and $\cT f\colon \KK\to\RR$ by
\begin{equation}\label{E-muT}
\mu \cT(A) \,\df\, \int_{\KK} \mu(\D{x},\D{u}) \cT( A \smid x, u)\,,
\quad A\in\Bor(\XX)\,,
\end{equation}
and
\begin{equation}\label{E-cTf}
\cT f(x,u) \,\df\, \int_{\XX} f(y) \cT( \D{y} \smid x, u)\,,
\quad (x,u)\in\KK\,,
\end{equation}
respectively. We use the convenient notation for integrals of functions
\begin{equation}\label{E-not1}
\mu(f)\,=\,\langle \mu, f\rangle \,\df\, \int_{\KK} f(x,u)\,\mu(\D{x},\D{u})\,,
\end{equation}
and similarly for $f\in\cM_b(\XX)$ and $\mu\in\Pm(\XX)$ if no ambiguity arises. We then have
\begin{equation*}
\langle \mu\cT, f\rangle\,=\, \langle \mu, \cT f\rangle
\qquad\text{for\ } \mu\in\Pm(\KK)\,,\ f\in\cM_b(\XX)\,.
\end{equation*}

The set of \emph{invariant occupation measures} (or, as is used more commonly in the literature: \emph{ergodic occupation measures} \cite{survey}) is defined by
\begin{equation}\label{E-cG}
\cG \,\df\,
\bigl\{\mu\in\Pm(\KK)
\colon \mu(B \times\Act) = \mu \cT( B), \ B \in \Bor(\XX) \bigr\}\,.
\end{equation}
We also let
\begin{equation*}
\cH \,\df\,
\biggl\{\uppi \in \Pm(\XX)
\colon \exists \gamma \in \Gamma_{S} \text{\ such that\ } \uppi(A) = \int_\XX
\cT^{\gamma}(A\smid x)\,\uppi(\D{x}),\  A \in \Bor(\XX) \biggr\}
\end{equation*}
denote the set of \emph{invariant probability measures} of the controlled
Markov chain.

Let $\mu\in\cG$.  It is well known that $\mu$ can be disintegrated into
a stochastic kernel $\phi$ on $\XX\times\Bor(\Act)$ and $\uppi\in\Pm(\XX)$ such that
\begin{equation*}
\mu(\D{x},\D{u}) \,=\, \phi(\D{u}\smid x)\,\uppi(\D{x})\,,
\end{equation*}
and $\phi$ is $\uppi$-a.e. uniquely defined on the support of $\uppi$. We denote this disintegration by $\mu =(\phi \uppi)$.
Therefore, if $\gamma\in\Usm$ is any policy which agrees $\uppi$-a.e. with $\phi$,
then we have $\uppi(A) = \cT^{\gamma}(A\smid x)\,\uppi(\D{x})$
for $A\in\Bor(\XX)$.
Therefore, $\uppi\in\cH$.
Conversely, if $\uppi\in\cH$ with an associated $\gamma\in\Usm$, then it is clear from
the definitions that $(\gamma \uppi) \in\cG$.

\subsection{Two Complementary Regularity Assumptions}

\begin{itemize}
\item[\hypertarget{H1}{\textbf{(H1)}}]
The transition kernel $\cT$ is called \emph{weakly continuous} if the map
\begin{equation*}
\KK\,\ni\,(x,u)\mapsto\int_{\XX} f(z)\cT(\D{z}\smid x,u)
\end{equation*}
is continuous for all $f\in C_b(\XX)$.
\end{itemize}

\begin{itemize}
\item[\hypertarget{H2}{\textbf{(H2)}}]
The transition kernel $\cT$ satisfies the following:
\begin{itemize}
\item[(a)]
For any $x\in\XX$, the map $u\mapsto\int f(z)\cT(\D{z}\smid x,u)$ is continuous
for every bounded measurable function $f$.
\item[(b)] 
There exists a finite measure $\nu$ majorizing $\cT$, that is
\begin{equation}\label{EH2A}
\cT(\D{y} \smid x, u) \,\le\, \nu(\D{y})\,, \qquad x \in \XX,\; u \in \Act\,.
\end{equation}
\end{itemize}
\end{itemize}

\begin{example}\label{ornek11}
Consider the following model. 
\begin{equation*}
x_{n+1}\,=\, F(x_{n},u_{n})+w_{n}\,,\qquad n=0,1,2,\dotsc\,,
\end{equation*}
where $\mathbb{X}=\mathbb{R}^{n}$ and the $w_{n}$'s are independent and
identically distributed (i.i.d.) random vectors.

\begin{itemize}
\item[(i)] If $F$ is continuous, regardless of the distribution of $w_n$, \textbf{(H1)} applies.
\item[(ii)] If $u\mapsto F(x,u)$ is continuous for all $x \in \mathbb{X}$ and $w_n$ has a distribution which admits a bounded and continuous density function, then \textbf{(H2)}(a) holds. 
\item[(iii)] If in addition to (ii), we assume that $F$ is bounded, (\ref{EH2A}) applies.
\end{itemize}
\end{example}

When one associates with a controlled Markov model an average cost criterion with bounded cost functions, it has been shown that under these assumptions, for average cost criteria, optimal solutions exist under additional conditions (such as positive Harris recurrence and continuity in both the state and actions under \textbf{(H1)} \cite{survey,Borkar2,kurano1989existence,hernandez1993existence,survey,HernandezLermaMCP}, and continuity of the cost in actions for each state under \textbf{(H2)} \cite{arapostathis2021optimality}); generalizations to the unbounded setup follow from tightness/lower semi-continuity arguments.

%

\section{Young and Borkar Topologies on Stationary Control Policies}

Consider the deterministic measurable policy space (\ref{detPolS}) and the relaxed/randomized measurable policy space (\ref{ranPolS}). In the following, we will study two topologies on (\ref{ranPolS}).

\subsection{Young Topology on Control Policies}

As noted earlier, a common approach which has been ubiquitously adopted in various fields often with different terminologies, is via what is commonly defined as the Young topology approach (e.g., in optimal deterministic control \cite{young1937generalized,mcshane1967relaxed},\cite[Section 2.1]{castaing2004young},\cite[p. 254]{warga2014optimal},\cite{mascolo1989relaxation}; distributional strategies in economics \cite{milgrom1985distributional} \cite{mertens2015repeated}; and optimal quantization \cite{YukselOptimizationofChannels}). 

To present the Young topology on control policies, we first present a relevant representation result (see Borkar \cite{BorkarRealization}). 

Let $\mathbb{X}, \mathbb{M}$ be Borel spaces. Let ${\cal P}(\mathbb{X})$ denote the set of probability measures on $\mathbb{X}$. Consider the set of probability measures
\begin{eqnarray}\label{extremePointQuan0}
&&\Theta: = \bigg\{\zeta \in {\cal P}(\mathbb{X} \times \mathbb{M}): \nonumber \\
&& \qquad \quad \zeta(dx,dm) = P(dx) \, Q^f(dm|x), Q^f(\cdot | x) = 1_{\{f(x) \in \cdot\}}, f : \mathbb{X} \to \mathbb{M} \bigg\}
\end{eqnarray}
on $\mathbb{X} \times \mathbb{M}$ with fixed input marginal $P$ on $\mathbb{X}$ and with the stochastic kernel from $\mathbb{X}$ to $\mathbb{M}$ realized by any measurable function $f : \mathbb{X} \to \mathbb{M}$. We equip this set with the weak convergence topology. This set is the (Borel measurable) set of the extreme points of the set of probability measures on $\mathbb{X} \times \mathbb{M}$ with a fixed marginal $P$ on $\mathbb{X}$. For compact $\mathbb{M}$, the Borel measurability of $\Theta$ follows \cite{Choquet} since the set of probability measures on $\mathbb{X} \times \mathbb{M}$ with a fixed marginal $P$ on $\mathbb{X}$ is a {\it convex and compact} set in a complete separable metric space, and therefore, the set of its extreme points is Borel measurable; measurability for the non-compact case follows from \cite[Lemma 2.3]{BorkarRealization}.
Furthermore, given a fixed marginal $P$ on $\mathbb{X}$, any stochastic kernel $Q$ from $\mathbb{X}$ to $\mathbb{M}$ can almost surely be identified by a probability measure $\Xi \in {\cal P}(\Theta)$ such that
\begin{eqnarray}\label{convR0}
Q(\cdot|x)  = \int_{\Theta} \Xi(dQ^f) \, Q^f(\cdot|x). 
\end{eqnarray}
In particular, a stochastic kernel can thus be viewed as an integral representation over probability measures induced by deterministic policies (i.e., a mixture of deterministic policies). 


Let $\mathbb{X}$ be a Polish space and let $\mathcal{P}(\mathbb{X})$ denote the family of
all probability measures on $(\mathbb{X},\mathcal{B}(\mathbb{X}))$. Let $\{\mu_n,\, n\in \mathbb{N}\}$ be a sequence in
$\mathcal{P}(\mathbb{X})$.

The sequence $\{\mu_n\}$ is said to  converge
  to $\mu\in \mathcal{P}(\mathbb{X})$ \emph{weakly} if
\begin{align}\label{weakConvD}
 \int_{\mathbb{X}} c(x) \mu_n(dx)  \to \int_{\mathbb{X}}c(x) \mu(dx)
\end{align}
for every continuous and bounded $c: \mathbb{X} \to \mathbb{R}$.

Let $\mathbb{Y}$ be another Polish space.
\begin{definition}\label{swTopDefinition}
The $w$-$s$ (weak-setwise or setwise-weak) topology on the set of probability measures ${\cal P}(\mathbb{X} \times \mathbb{Y})$
is the coarsest topology under which
$\int f(x,y) \mu(dx,dy): {\cal P}(\mathbb{X} \times \mathbb{Y}) \to \mathbb{R}$
is continuous for every measurable and bounded $f$ which is continuous in $y$
for every $x$ (but unlike the weak topology, $f$ does not need to be continuous in $x$). Let $\{\mu_n,\, n\in \mathbb{N}\}$ be a sequence in $\mathcal{P}(\mathbb{X} \times \mathbb{Y})$. The sequence $\{\mu_n\}$ is said to  converge
  to $\mu\in {\cal P}(\mathbb{X} \times \mathbb{Y})$ in the $w$-$s$ topology if
  \[\int f(x,y) \mu_n(dx,dy) \to \int f(x,y) \mu(dx,dy),\]
  for every measurable and bounded $f$ which is continuous in $y$ for every $x \in \mathbb{X}$.
\end{definition}
See \cite[Theorem 3.7]{schal1975dynamic} for some properties and equivalence relations on the $w$-$s$ convergence. Note that the coordinate on which continuity is imposed should be made explicit (thus, the $w$-$s$ or $s$-$w$ topologies are identically denoted).

\sym{\begin{definition}{\bf Convergence of Policies in $\Gamma_S$ under Young topology at reference (input) measure $\mu$.}  Let $\mu$ be a $\sigma$-finite measure. A sequence of stationary policies $\gamma_n \to \gamma \in \Gamma_S$ under Young topology at input $\mu$, if the joint measure $(\mu \gamma_n) \to (\mu\gamma)$ in the $w$-$s$ sense at input $\mu$, i.e., for every measurable and bounded $g: \mathbb{X} \times \mathbb{U} \to \mathbb{R}$ with $g(x, \cdot): u \mapsto g(x,u)$ continuous and
\[\int \mu(dx) \sup_{u \in \mathbb{U}} |g(x,u)| < \infty,\]
\begin{align}\label{convWeakSet}
\int \mu(dx) \bigg(\gamma_n(du|x) g(x,u) \bigg) \to \int \mu(dx) \bigg( \gamma(du|x) g(x,u) \bigg)  
\end{align}

\end{definition}
}


\begin{lemma}\label{wsYoungCon}
If the probability measures induced by the two-tuple random variable sequence $(X,U_n)$ converges to that induced by $(X,U)$ weakly, then this convergence is also in the $w$-$s$ (setwise-weak) sense (setwise in $x$, weakly in $u$).
\end{lemma}
\textbf{Proof.}
The marginal on $X$ is fixed along the sequence. The result then follows from \cite[Theorem 3.7]{schal1975dynamic} (or \cite[Theorem 2.5]{balder2001}).
\qed

\begin{remark}
\sym{In view of Lemma \ref{wsYoungCon}, when $\mu$ is a probability (or more generally, a finite) measure, since the marginal of the joint measure $(\mu \gamma_n)$ on $\mathbb{X}$ is fixed, the convergence in (\ref{convWeakSet}) is implied by weak convergence as well (with $g(\cdot,\cdot)$ taken to be bounded continuous on $\mathbb{X} \times \mathbb{U}$). We remark also that the above equivalence does not necessarily apply when the input $\mu$ is not a finite measure.} 
\end{remark}

In the literature, further related results are present; see \cite{balder1997consequences} for a concise discussion and \cite[Theorem 2.2]{balder1988generalized} and \cite[Theorem 3.32]{florescu2012young} for equivalent characterizations of Young topology convergence.

\subsection{Borkar Topology on Control Policies}

We now present a different topology, introduced by Borkar \cite{borkar1989topology} \cite[Lemma 2.4.1]{arapostathis2012ergodic}. 

As we noted earlier, the presentation in \cite[Section 3]{borkar1989topology} and \cite[Section 2.4]{arapostathis2012ergodic} are slightly different though the induced topologies are equivalent. Saldi \cite{SaldiArXiv2017} generalized Borkar's topology to setups where $\mathbb{X}$ is a general standard Borel space. In the following, we will take $\mathbb{U}$ to be compact, though the results generalize when $\mathbb{U}$ is locally compact. 

Consider the relaxed control policy space (\ref{ranPolS}). The Borkar topology (and its generalization by Saldi) on this space is formulated as the weak$^*$ topology when (\ref{ranPolS}) is viewed as maps from $\mathbb{X}$ to the space of signed measures ${\cal M}(\mathbb{U})$ with finite variation of which probability measures ${\cal P}(\mathbb{U})$ is a subset, and signed measures with finite variation are studied as the topological duals of continuous functions vanishing at infinity (see \cite[Theorem 7.17]{Fol99} and \cite[Theorem 1.5.5, p. 27]{CeMe97}). By the Banach-Alaoglu theorem \cite[Theorem 5.18]{Fol99}, as the unit ball of $L_{\infty}(\mathbb{R}, {\cal M}(\mathbb{U}))$ is compact under the weak$^*$ topology, this leads to a compact metric topology on relaxed control policies. An equivalent representation of this topology is given by the following \cite[Lemma 2.4.1]{arapostathis2012ergodic}:

\begin{definition}\label{BorkarTopDef} \cite{borkar1989topology} \cite[Lemma 2.4.1]{arapostathis2012ergodic}
{\bf Convergence of Policies in $\Gamma_S$ under Borkar topology.} With $\mathbb{X}=\mathbb{R}^d$, a sequence of stationary policies $\gamma_n \to \gamma \in \Gamma_S$ in the Borkar topology if for every continuous and bounded $g: \mathbb{X} \times \mathbb{U} \to \mathbb{R}$ and every $f \in L_1(\mathbb{X}) \cap L_2(\mathbb{X})$
\begin{align}\label{convWeakSet11}
\int f(x) \int \gamma_n(du|x) g(x,u) dx \to \int f(x) \int \gamma(du|x) g(x,u) dx 
\end{align}
\end{definition}

\begin{remark}
Building on, and slightly revising, the analysis in the proof of Lemma \ref{wsYoungCon}, the functions $g$ in Definition \ref{BorkarTopDef} may be relaxed to be continuous only in $u$ for every $x \in \mathbb{X}$. 
\end{remark}

The generalization of this topology by Saldi \cite{SaldiArXiv2017} is the following, where the input space $\mathbb{X}$ is arbitrary standard Borel, though with a fixed input measure $\mu$. Let $C_0(\mathbb{U})$ be the Banach space of all continuous real functions on $\mathbb{U}$ vanishing at infinity endowed with the norm
\begin{align}
\|g\|_{\infty} = \sup_{u \in \mathbb{U}} |g(u)|. \nonumber
\end{align}
Saldi formulated this topology via noting that with $L_1\bigl(\mu,C_0(\mathbb{U})\bigr)$ denoting the set of all Bochner-integrable functions from $\mathbb{X}$ to $C_0(\mathbb{U})$ endowed with the norm
\begin{align}
\|f\|_1 \coloneqq \int_{\mathbb{X}} \|f(x)\|_{\infty} \, \mu(dx), \nonumber
\end{align}
using the fact that $C_0(\mathbb{U})^* = {\cal M}(\mathbb{U})$, and that the topological dual of $\left(L_1\bigl(\mu,C_0(\mathbb{U}\bigr),\|\cdot\|_1\right)$ can be identified with $\left(L_{\infty}\bigl(\mu,{\cal M}(\mathbb{U})\bigr),\|\cdot\|_{\infty}\right)$ \cite[Theorem 1.5.5, p. 27]{CeMe97} (see also \cite{dieudonne1951theoreme,fattorini1994existence} for further context on such duality results, in particular when instead of countably additive signed measures, finitely additive such measures are considered); that is,
$$
L_1\bigl(\mu,C_0(\mathbb{U})\bigr)^* = L_{\infty}\bigl(\mu,{\cal M}(\mathbb{U})\bigr).
$$

\begin{definition} \cite{SaldiArXiv2017}
{\bf Convergence of Policies in $\Gamma_S$ under Borkar-Saldi topology at input $\mu$.} With $\mathbb{X}$ standard Borel, a sequence of stationary policies $\gamma_n \to \gamma \in \Gamma_S$ in the Borkar-Saldi topology if for every measurable bounded $g: \mathbb{X} \times \mathbb{U} \to \mathbb{R}$ with $g(x , \cdot): \mathbb{U} \to \mathbb{R}$ continuous and 
\[\int \mu(dx) \sup_{u \in \mathbb{U}} |g(x,u)| < \infty,\]
we have
\begin{align}\label{convWeakSet1111}
\int \mu(dx) \int \gamma_n(du|x) g(x,u) \to \int \mu(dx) \int \gamma(du|x) g(x,u)  
\end{align}
\end{definition}

Note that this definition is closely related with Young topology convergence at input $\mu$, and when $\mathbb{U}$ is compact this can be seen to be identical to the Young topology convergence at input $\mu$ (via a uniform integrability argument following the steps to be given in the proof of Lemma \ref{equivYoungT}). A particular difference with the Borkar topology, and similarity with Young topology, is with regard to the presence of a fixed input measure.

The Borkar topology has been used to show that, for non-degenerate controlled diffusions, expected cost is continuous in control policies under this topology (see \cite{borkar1989topology} for finite horizon problems and \cite{pradhan2022near} for infinite horizon problems as well as continuity in Markov policies in addition to stationary policies) and for the continuity of invariant measures for such diffusions in control policies \cite[Lemma 4.4]{arapostathis2010uniform}. The generalization by Saldi was utilized to study existence of team optimal policies in decentralized stochastic control \cite{SaldiArXiv2017}.


\subsection{Some Properties and Equivalence Conditions of Young and Borkar topologies}

We first present a supporting result on convergence under the Young topology. 
\begin{lemma}\label{equivYoungT}
Let $\eta \ll \kappa$, where $\kappa$ is $\sigma$-finite and $\eta$ is a finite measure. Then, $\gamma_n \to \gamma$ under Young topology at input $\kappa$ implies $\gamma_n \to \gamma$ (under Young topology) at input $\eta$.
\end{lemma}
\textbf{Proof.}
By assumption we have that there exists $f: \mathbb{X} \to \mathbb{R}_+$ with
\[f(x) = \frac{d \eta}{d \kappa}(x),\]
which is the Radon-Nikodym derivative. In particular, $f$ is integrable under $\kappa$. 

Now, let $\bar{g}$ be continuous in $u$ for every $x$ and bounded where $\sup_{u \in \mathbb{U}}|\bar{g}(x,u)|$ is integrable under $\kappa$. Consider a sequence of converging control policies at $\kappa$ so that:
\[\lim_{n \to \infty} \int (\kappa\gamma_n)(dx,du) \bar{g}(x,u) = \int \kappa\gamma(dx,du) \bar{g}(x,u) \]

We will show that the above implies that for every continuous and bounded $g$ (which is necessarily integrable under the finite measure $\eta$)
\[\lim_{n \to \infty} \int (\eta\gamma_n)(dx,du) g(x,u) = \int \eta\gamma(dx,du) g(x,u) \]

We have the following uniform integrability condition:
\begin{align}
& \lim_{M \to \infty} \sup_{n \in \mathbb{N}} \int_{\{|f(x)| \geq M\}} \kappa(dx) \bigg|f(x) \bigg( \int \gamma_n(du|x)  g(x,u) \bigg) \bigg| \nonumber \\ 
&\leq \lim_{M \to \infty} \sup_{n \in \mathbb{N}} \int_{\{|f(x)| \geq M\}} \kappa(dx) f(x) \|g\|_{\infty} = 0.\label{unifInt}
\end{align}
Therefore, with
\begin{align}
&\int (\kappa\gamma_n)(dx,du) f(x) g(x,u) \nonumber  \\
&  = \int (\kappa\gamma_n)(dx,du) (f(x)- \min(M,f(x))) g(x,u)  \nonumber \\
& \quad \quad \quad + \int (\kappa\gamma_n)(dx,du) \min(M,f(x)) g(x,u)
\end{align}
the first term on the right is uniformly and absolutely bounded by $\epsilon_M$ with $\epsilon_M \to 0$ as $M \to \infty$ (by uniform integrability in (\ref{unifInt})), and the second term on the right converges to $\int \kappa\gamma(dx,du) \min(M,f(x)) g(x,u)$, for every $M \in \mathbb{R}_+$ (since for every $M$, $\min(M,f(x)) g(x,u)$ is bounded and continuous in $u$ for every $x$). Therefore,
\begin{align}
&\limsup_{n \to \infty} \int (\kappa\gamma_n)(dx,du) f(x) g(x,u) \nonumber  \\
&  \leq \epsilon_M + \int (\kappa\gamma) (dx,du) \min(M,f(x)) g(x,u) \leq  \epsilon_M + \int (\kappa\gamma)(dx,du) f(x) g(x,u)
\end{align}
and
\begin{align}
&\liminf_{n \to \infty} \int (\kappa\gamma_n)(dx,du) f(x) g(x,u) \nonumber  \\
&  \geq - \epsilon_M + \int (\kappa\gamma) (dx,du) \min(M,f(x)) g(x,u)
\end{align}
Since the above holds for every $M$, taking $M \to \infty$ and applying the dominated convergence theorem, the limit infimum and supremum of the sequences are seen to be equal. We arrive at
\begin{align}
& \int \eta \gamma_n(dx,du)  g(x,u)  = \int (\kappa\gamma_n)(dx,du) f(x) g(x,u)  \nonumber \\
& \qquad \qquad \to \int (\kappa\gamma)(dx,du) f(x) g(x,u) = \int (\eta \gamma)(dx,du)  g(x,u) 
\end{align}
\qed

\begin{remark}
Lemma \ref{equivYoungT} can be generalized to a uniform convergence over collection of finite input measures $\Upsilon=\{\eta_{\alpha}\}$ where $\eta_{\alpha} \ll \kappa$ with uniformly integrable (under $\kappa$) Radon-Nikodym derivatives $\frac{d\eta_{\alpha}}{d\kappa}$. In particular, a close look at the proof reveals that we would have, for all bounded continuous $g$:
\[\lim_{n \to \infty} \sup_{\eta_{\alpha} \in \Upsilon} \bigg| \int \eta_{\alpha} \gamma_n(dx,du)  g(x,u) - \int \eta_{\alpha} \gamma(dx,du)  g(x,u) \bigg| = 0\]
This has implications for convergence of policies when the underlying input measure is varying.  
\end{remark}


An implication of Lemma \ref{equivYoungT}, essentially a corollary, is then the following.

\begin{theorem}\label{equivYouBorTop}
 Let $\mathbb{X} = \mathbb{R}^n$ and $\lambda$ be the Lebesgue measure. Consider convergence in Young topology at some $\sigma$-finite input measure $\psi$.
 \begin{itemize}
 \item[(i)] If $\psi \ll \lambda$ with $h(x)=\frac{d\psi}{d\lambda}(x)$ is positive everywhere, then convergence in Young topology at input measure $\psi$ implies convergence in Borkar topology.
 
 \item[(ii)] If $\psi \ll \lambda$ and $\psi$ is a finite measure, then convergence in Borkar topology implies convergence in Young topology at input $\psi$.
  
 \item[(iii)] Convergence in Borkar topology is equivalent to convergence in Young topology at input $\psi$ if $\psi \ll \lambda$, $\psi$ is finite, and $h(x)=\frac{d\psi}{d\lambda}(x)$ is positive everywhere.
 \end{itemize}

\end{theorem} 

\textbf{Proof.} (i) Since $L_1$ is separable, the convergence in the Borkar topology is equivalent to convergence of
\begin{align}\label{convWeakSet112}
\int f_m(x) \int \gamma_n(du|x) g(x,u) \lambda(dx) \to \int f_m(x) \int  \gamma(du|x) g(x,u) \lambda(dx) 
\end{align}
for a countable collection of functions $\{f_m, m=1,2,\cdots\}$. Equivalently,
\begin{align}\label{convWeakSet12}
\int \frac{f_m(x)}{\|f_m\|_1} \int \gamma_n(du|x) g(x,u) \lambda(dx) \to \int  \frac{f_m(x)}{\|f_m\|_1} \int   \gamma(du|x) g(x,u) \lambda(dx) 
\end{align}

For each fixed $m$, the convergence of (\ref{convWeakSet12}) for every bounded and continuous $g$ can be viewed as weak convergence of the signed measure $(F_m \gamma_n)(A,B) = \int_A \frac{f_m(x)}{\|f_m\|_1} \int_B \gamma_n(du|x) dx$, where we define this space as a locally convex space with the semi-norms defined by convergences of
\[\bigg|\int (F_m \gamma_n)(dx,du) g(x,u) - \int (F_m \gamma)(dx,du) g(x,u) \bigg|\]
to zero for each continuous and bounded $g$.

Now, define another signed measure ${\bf F}$ by:
\[{\bf F}(A) = \sum_m 2^{-m} \int_A \frac{f_m(x)}{\|f_m\|_1} dx\] 
This is a finite measure and one that satisfies $F_m \ll {\bf F}$. 

Thus, by Lemma \ref{equivYoungT} it follows that Young convergence at ${\bf F}$ will imply the Borkar convergence topology. 

Let $\psi \ll \lambda$ with $h(x)=\frac{d\psi}{d\lambda}(x)$ be positive everywhere. It follows then that ${\bf F} \ll \psi$ and Young topology convergence at $\psi$ implies convergence in the Borkar topology.  

(ii) For any given finite measure ${\bf F}$ which is absolutely continuous with respect to the Lebesgue measure, the Radon-Nikodym derivative $\frac{d{\bf F}}{d\lambda}$ (which is necessarily in $L_1(\mathbb{X})$) can be approximated arbitrarily well by functions $f_m$ in $L_1(\mathbb{X}) \cap L_2(\mathbb{X})$, in the $L_1$ norm. Therefore, by (\ref{convWeakSet112}), convergence in the Borkar topology implies that under the Young topology in this case.

(iii) This follows from (i) and (ii).

\qed

\begin{remark}
If the reference measure $\psi$ is not absolutely continuous with respect to the Lebesgue measure, then it is not necessarily the case that the Borkar topology implies Young's: Any atomic probability measure can be approximated by nonatomic measures under the weak convergence topology, but this is not strong enough for the convergence in Young topology as control policies would also be varying along a sequence of policies (unless one restricts the policies to be continuously converging (see \cite[Theorem 3.5]{serfozo1982convergence} or \cite[Theorem 3.5]{Lan81}).
\end{remark}

\section{Continuous Dependence of Invariant Measures on Stationary Control Policies under Young Topology}

In this section, we present the following continuity results involving invariant measures on the space of control policies. We first consider the case with $\mathbb{X} = \mathbb{R}^n$. 
\begin{theorem}\label{ContInvPolicy}
Suppose that 
\begin{itemize}
\item[\ttup{i}]
We have $\mathbb{X}=\mathbb{R}^n$ for some finite $n$, and for all $x \in \mathbb{R}^n$, $\mathbb{U}(x) = \mathbb{U}$ is compact.
\item[\ttup{ii}]
For every stationary policy $\gamma \in \Gamma_{\mathsf S}$ there exists a unique invariant probability measure. 

%
\item[\ttup{iii}] The kernel $\cT(\D{y} \smid x,u)$ is such that, the family of conditional probability measures $\{\cT(\D{y} \smid x,u), x \in \mathbb{X}, u\in \mathbb{U}\}$ admit densities $f_{x,u}(y)$ with respect to a reference measure $\psi$, and all such densities are bounded and equicontinuous (over $x \in \mathbb{X}, u\in \mathbb{U}$). 

\item[\ttup{iv}] One of the following holds: (a) \hyperlink{H1}{(H1)} holds and $\cG$ (defined in (\ref{E-cG})) is weakly compact, or (b) \hyperlink{H2}{(H2)} holds (for which $\cG$ is necessarily weakly compact \cite[Theorem 2.2]{arapostathis2021optimality}). 

\end{itemize}
Then, the invariant measure is weakly (and in total variation) continuous on the control policy space $\Gamma_S$, where we endow the stationary policy space $\Gamma_S$ with the Young topology at reference (input) measure $\psi$. 
\end{theorem}


\textbf{Proof.}


Observe that the family of densities $f_{x,u}(\cdot)$ being equicontinuous over $x \in \mathbb{X}, u\in \mathbb{U}$ implies that $\{ \int f_{x,u}(y) \mu(dx)\gamma(du|x), \quad \mu \in {\cal P}(\mathbb{X}), \gamma \in \Gamma_{\mathsf S}\}$ is also equicontinuous. Thus, the family of invariant probability measures under any stationary policy admit densities, with respect to the reference probability measure $\psi$, which are bounded and equicontinuous.
Then, following e.g., \cite[Lemma 4.3]{YukselOptimizationofChannels}, if $\{h_n\}$ is a sequence of probability density functions (with respect to some reference measure $\psi$) which are equicontinuous and uniformly bounded and if the corresponding sequence of measures $\mu_n(dy) = h_n(y)\psi(dy) \to \mu(dy) = h(y) \psi(dy)$ weakly, then as a consequence of the Arzel\'a-Ascoli theorem (applied to $\sigma$-compact spaces) $h_n \to h$ pointwise and by Scheff\'e’s theorem \cite{Bil86}, $\mu_n \to \mu$ in total variation. 

Let $\gamma$ be any randomized stationary policy. Suppose that this policy gives rise to an invariant probability measure $\uppi_\gamma(\D{x},\D{u})$ (by (ii)). Now, consider a sequence of policies $f_n$ so that under this sequence of policies $f_n \to \gamma$ at $\psi$ and therefore, by Lemma \ref{equivYoungT}, $\uppi_\gamma(\D{x}) f_n(\D{u}|x)$ converges weakly to $\uppi_\gamma(\D{x},\D{u})$ as well.

Now, let us apply the same control policy sequence to the random variable $X_n$ which has the probability measure $\uppi_{f_n}(\D{x})$ equal to the marginal of the invariant measure under policy $U=f_n(X)$. Then, for every continuous and bounded $g \in C_b(\XX)$
\begin{equation}\label{denk1Inv}
\int \uppi_{f_{n}}(\D{x})f_n(\D{u}|x) \bigg( \int g(y) \cT(\D{y} \smid  x,u) \bigg) = \int \uppi_{f_{n}}(\D{x}) g(x)\,.
\end{equation}
Let $\uppi_{f_{n_k}}(\D{x})$ be a weakly converging subsequence with limit $\eta$ (by the compactness assumption on $\cG$). By hypothesis (as a result of the discussion above on Scheff\'e’s theorem), this convergence is also in total variation. Therefore, via writing (by dropping the subsequence notation $n_k$ and considering this as a converging sequence)
\begin{align}\label{denk2Inv}
&\int \uppi_{f_{n}}(\D{x})f_n(\D{u}|x)  \bigg( \int g(y) \cT(\D{y} \smid  x,u) \bigg) -\int  \eta(\D{x})\gamma(\D{u}|x) \bigg( \int g(y) \cT(\D{y} \smid  x,u) \bigg)  \nonumber \\
& = \int \uppi_{f_{n}}(\D{x})f_n(\D{u}|x)  \bigg( \int g(y) \cT(\D{y} \smid  x,u) \bigg) - \int  \eta(\D{x})f_n(\D{u}|x)  \bigg( \int g(y) \cT(\D{y} \smid  x,u) \bigg) \nonumber \\
& \quad + \int \eta(\D{x})f_n(\D{u}|x) \bigg( \int g(y) \cT(\D{y} \smid  x,u) \bigg) -  \int \eta(\D{x})\gamma(\D{u}|x)  \bigg( \int g(y) \cT(\D{y} \smid  x,u) \bigg),
\end{align} 
we have that for every continuous and bounded function $g$ (using the inner-product $\langle \cdot, \cdot \rangle$ notation for the above integrations)
\begin{align}\label{ConvProofDenk11}
 \bigg\langle \uppi_{f_{n}}(\D{x})f_n(\D{u}|x) - \eta(\D{x})f_n(\D{u}|x) , \bigg( \int g(y) \cT(\D{y} \smid  x,u) \bigg)  \bigg\rangle \to 0
 \end{align}
and
\begin{align}\label{ConvProofDenk12}
\bigg\langle \eta(\D{x})f_n(\D{u}|x) - \eta(\D{x})\gamma(\D{u}|x), \bigg( \int g(y) \cT(\D{y} \smid  x,u) \bigg)  \bigg\rangle \to 0
\end{align}

The first term (\ref{ConvProofDenk11}) here on the right converges to zero due to total variation convergence of $\uppi_{f_{n}}$ to $\eta$ (since we apply the same policy $f_n$, and convergence is uniform over all measurable functions as in the proof of \cite[Lemma 1.1(iii)]{KYSICONPrior}). 

We now show that the second term (\ref{ConvProofDenk12}) converges to zero as well. By the Arzel\'a-Ascoli theorem implying that the limit function defines a density with respect to $\psi$, we have $\eta \ll \psi$. Therefore, by the definition of convergence (of control policies at $\psi$, the reference measure for each invariant probability measure under any stationary  policy) and Lemma \ref{equivYoungT}, we have that the second term converges to zero: Recall that
\[\bigg( \int g(y) \cT(\D{y} \smid  x,u) \bigg)\]
is either continuous in $x,u$ (under  \hyperlink{H1}{(H1)}); or continuous in $u$ for every $x \in \mathbb{X}$ (under \hyperlink{H2}{(H2)(a)}). In either case, by \cite[Theorem 3.10]{schal1975dynamic} or \cite[Theorem 2.5]{balder2001}, since the measure converges weakly and the marginal on $\mathbb{X}$ converges setwise, the convergence is also in the $w$-$s$ sense (see Lemma \ref{wsYoungCon}). 


The term on the right of (\ref{denk1Inv}) converges to $\int \eta(\D{x}) g(x)$, due to convergence of the measures in total variation. Thus, we have that the following holds for every continuous and bounded $g$:
\begin{equation*}
\int \eta(\D{x})\gamma(\D{u}|x) \bigg( \int g(y) \cT(\D{y} \smid  x,u) \bigg) = \int \eta(\D{x}) g(x)\,.
\end{equation*}
Since continuous and bounded functions are measure determining (see \cite[p. 13]{Billingsley} or \cite[Theorem 3.4.5]{ethier2009markov}), the above implies that $\eta$ is invariant.

However, by uniqueness of the invariant probability measure for any stationary policy, we have that $\eta$ must be $\uppi_{\gamma}$.

\qed

We remark that compactness of invariant measures under total variation has been studied in \cite[Lemma 3.2]{borkarghosh1990controlled} with a similar argument in continuous-time.

Let us present some concrete examples. Assume either Example \ref{ornek11}(i) or Example \ref{ornek11}(ii) holds, together with Example \ref{ornek11}(iii), and assume that $w_n$ admits a bounded density function which is positive everywhere (so that the induced Markov chain by any stationary policy is irreducible) and has a bounded derivative (so that the densities are equicontinuous); a Gaussian density is sufficient for each of these conditions. In this case, Theorem \ref{ContInvPolicy} is applicable.

The following result, essentially providing additional examples for which Theorem \ref{ContInvPolicy} is applicable, is a corollary of the above, since with finite spaces the continuity conditions (weak or setwise) are equivalent, and every probability measure is absolutely continuous with respect to a weighted counting measure; and thus the proof follows identically.

\begin{corollary}\label{ContInvPolicy2}
\begin{itemize}
\item[(i)] Let $\mathbb{X}, \mathbb{U}$ be finite and for every stationary policy there be a unique invariant probability measure. Then, the invariant measure is continuous under Young topology on stationary controls at any input which places positive mass on each element in $\mathbb{X}$. 
\item[(ii)] The above applies to the case where $\mathbb{X}$ is countable and $\mathbb{U}$ is compact, provided that ${\cal G}$ (\ref{E-cG}) is compact and the transition kernel is continuous in actions.
\end{itemize}
\end{corollary}

\begin{remark}
With finite or countable spaces, equivalent to the Young topology with an input measure placing positive mass on each state $x \in \mathbb{X}$, one can also directly work with distance measures of the form:
\[ \sum_{x \in \mathbb{X}} 2^{-x} \rho(\gamma_n(du|x), \gamma(du|x)),\]
with $\rho$ being any weak convergence inducing metric. Such distance measures have been shown to be consequential in stochastic dynamic game theory and learning theoretic algorithms (see e.g. \cite[Section 2.4]{yongacoglu2021satisficing}). 
\end{remark}

\section{An Application: Near Optimality of Continuous or Quantized (in state and action) Deterministic Policies in Average Cost Optimal Control}

Consider the following average cost problem of finding
\begin{eqnarray}\label{AverageCostProblemDef}
J^*(x):= \inf_{\gamma} J(x,\gamma) = \inf_{\gamma \in \Gamma_A} \limsup_{T \to \infty} {1 \over T} E^{\gamma}_x [\sum_{t=0}^{T-1} c(x_t,u_t)]
\end{eqnarray}

It has been shown that, for such average cost criteria, optimal solutions exist under additional conditions (such as positive Harris recurrence, boundedness of the cost function and its continuity in both the state and actions under \textbf{(H1)} \cite{survey,Borkar2,kurano1989existence,hernandez1993existence,survey,HernandezLermaMCP}, or continuity of the cost only in actions under \textbf{(H2)} \cite{arapostathis2021optimality}; and a tightness condition on the set of invariant probability measures induced by stationary policies).

We present an approximation result. 

\begin{theorem}\label{denseQuantizedFiniteAC}
Let the conditions of Theorem \ref{ContInvPolicy} hold, either with \hyperlink{H1}{(H1)} or \hyperlink{H2}{(H2)(a)}. Further assume that $\psi$ is a finite measure (a sufficient condition being \hyperlink{H2}{(H2)(b)}). Suppose that the measurable cost function is bounded and is continuous in the actions for each state.
Then,
\begin{itemize}
\item[(i)] an $\epsilon$-optimal solution exists among continuous stationary (possibly randomized) control policies
\item[(ii)] an $\epsilon$-optimal solution exists among quantized (possibly randomized) stationary control policies which have finite range and finite domain with uniform quantization (i.e., with both state and action space quantization) 
\item[(iii)] if $\psi \ll \lambda$ with $\lambda$ being the Lebesgue measure, an $\epsilon$-optimal solution exists among continuous stationary deterministic control policies
\item[(iv)] if $\psi \ll \lambda$ with $\lambda$ being the Lebesgue measure, an $\epsilon$-optimal solution exists among stationary and deterministic policies which have finite range and finite domain with uniform quantization (i.e., with both state and action space quantization). 
\end{itemize}
\end{theorem}

In the following, we prove this theorem. Recall the following.

\begin{theorem}\cite[Theorem 7.5.2]{Dud02}\label{LusinDudley}[Lusin's Theorem]
Let $(\mathbb{X}, T)$ be any topological space and $\mu$ a finite,
closed regular Borel measure on $\mathbb{X}$. Let $(\mathbb{S}, d)$ be a separable metric space and let $f$ be a Borel-measurable function from $\mathbb{X}$ into $\mathbb{S}$. Then for any $\epsilon > 0$ there is a closed set $F \subset \mathbb{X}$ such that $\mu(\mathbb{X}\setminus F) < \epsilon$ and the restriction of $f$ to $F$ is continuous.
\end{theorem}

Finally we recall Tietze's extension theorem, which will be used in conjunction with Lusin's theorem. This will be used to construct a continuous extension of the continuous function defined on $F$ in Theorem \ref{LusinDudley} to $\mathbb{X}$. Note that we consider a general space setup in the following; this is needed as we will study stochastic kernels as probability measure-valued maps.

\begin{theorem}\cite[Theorem 4.1]{dugundji}\label{tietze}[Tietze's extension theorem]
Let $\mathbb{X}$ be an arbitrary metric space, $A$ a closed subset of $\mathbb{X}$, $L$ a locally convex linear space, and $f: A \rightarrow L$ a continuous map.
Then there exists a continuous function $f_C: \mathbb{X} \rightarrow L$ such that $f_C(a) = f(a) \: \forall a \in A$. Furthermore, the image of $f_C$ satisfies $f_C(\mathbb{X}) \subset$ closed convex hull of $f(A)$.
\end{theorem} 

\textbf{Proof of Theorem \ref{denseQuantizedFiniteAC}.}

We will first prove (i)-(ii), followed by (iii)-(iv).

By Theorem \ref{ContInvPolicy} (see also \cite[Theorems 2.1 and 2.2]{arapostathis2021optimality}, under the complementary conditions of \hyperlink{H1}{(H1)} or \hyperlink{H2}{(H2)(a)}) an optimal invariant measure exists, which leads to an optimal policy for almost every initial condition (over a set of measure one) under its corresponding invariant measure. We note that if the induced Markov chain is positive Harris recurrent, then the optimality would hold for all initial conditions given the assumptions. Consider this, possibly randomized, control policy $\gamma$, and consider the measure on the product space with input measure $\psi$. Now, note that the policy, as a stochastic kernel (as a regular conditional probability), is defined with the following property: for every $x \in \mathbb{X}$, $\gamma(du|x)$ is $\mathcal{P}(\mathbb{U})$-valued and for every Borel set $B \subset \mathbb{U}$, 
$\gamma(u \in B| \cdot): \mathbb{X} \to \mathbb{R}$ is Borel-measurable.

{\bf Step 1.}
Following \cite[Proposition 7.26]{BeSh78}, the above property is equivalent to the policy $\gamma(u^i \in \cdot|x)$ being a Borel-measurable map from $\mathbb{X}$ to $\mathcal{P}(\mathbb{U})$ (which is endowed by the weak convergence topology). Because $\mathbb{U}$ is standard Borel, $\mathcal{P}(\mathbb{U})$ is a separable metric space, and can be defined by viewing the space of probability measures $\mathcal{P}(\mathbb{U})$ as a convex subset of a locally convex space \cite[Chapter 3]{rudin1991functional} of \sym{signed measures with finite total variation} defined on ${\cal B}(\mathbb{U})$, where we define the locally convex space of signed measures with the following notion of convergence: We say that $\nu_n \rightarrow \nu$ if $\int f(u) \nu_n(du) \rightarrow \int f(u) \nu(du)$ for every continuous and bounded function $f: \mathbb{U} \to \mathbb{R}$. 

{\bf Step 2.}
In the following, we will first approximate $\gamma$ with a continuous $\gamma_C$ which will then be approximated with policies which are quantized.

Theorems \ref{LusinDudley} and \ref{tietze} apply with the continuous extension being probability measure-valued (by Theorem \ref{tietze}). Accordingly, for every $\epsilon>0$, we have a continuous map $\gamma_C$ which agrees with $\gamma$ except on a set $K_{\epsilon}$ of measure $\psi(K_{\epsilon}) \leq \epsilon$.

Now, we find a sequence of continuous functions $\gamma^n_C$ which converges to $\gamma$ in Young topology at input $\psi$. To see this, note that convergence is equivalent to 
\[ \int_{\mathbb{X}} \int_{\mathbb{U}} \psi(dx) \gamma^n_C(du|x) g_m(x,u) \to \int_{\mathbb{X}} \int_{\mathbb{U}} \psi(dx) \gamma(du|x) g_m(x,u),\]
for a countable collection of continuous and (uniformly) bounded functions, as these are measure determining (see \cite[p. 13]{Billingsley} or \cite[Theorem 3.4.5]{ethier2009markov}). 

We can find, for every $\epsilon_n > 0$, a set $K_{\epsilon_n}$ so that on the complement of this set $\gamma^n_C$ agrees with $\gamma$. In particular, simultaneously for every $m$,
\begin{align}
&\bigg| \int_{\mathbb{X}} \int_{\mathbb{U}} \psi(dx) \gamma^n_C(du|x) g_m(x,u) - \int_{\mathbb{X}} \int_{\mathbb{U}} \psi(dx) \gamma(du|x) g_m(x,u) \bigg| \nonumber \\
&=\bigg| \int_{K_{\epsilon_n}} \int_{\mathbb{U}} \psi(dx)\gamma^n_C(du|x) g_m(x,u) - \int_{K_{\epsilon_n}} \int_{\mathbb{U}} \psi(dx) \gamma(du|x) g_m(x,u) \bigg| \nonumber \\
&\leq  \int_{K_{\epsilon_n}} \psi(dx)  \bigg| \int_{\mathbb{U}} \gamma^n_C(du|x) g_m(x,u) -  \int_{\mathbb{U}}  \gamma(du|x) g_m(x,u)\bigg| \nonumber \\
& \leq 2\|g_m\|_{\infty} \epsilon_n \leq 2 M \epsilon_n,
\end{align}
where $M = \sup_{m \in \mathbb{N}} \|g_m\|_{\infty}$.
We can define a metric to determine convergence (of such continuous functions $\gamma_C$ to $\gamma$), and obtain:
\[d(\gamma_C,\gamma) := \sum_{m \in \mathbb{N}} 2^{-m} \frac{ \bigg| \int_{\mathbb{X}} \int_{\mathbb{U}} \psi(dx) \gamma_C(du|x) g_m(x,u) - \int_{\mathbb{X}} \int_{\mathbb{U}} \psi(dx) \gamma(du|x) g_m(x,u) \bigg|}{1+\bigg| \int_{\mathbb{X}} \int_{\mathbb{U}} \psi(dx) \gamma_C(du|x) g_m(x,u) - \int_{\mathbb{X}} \int_{\mathbb{U}} \psi(dx) \gamma(du|x) g_m(x,u) \bigg|}, \]
so that
\[d(\gamma^n_C,\gamma) \leq M \epsilon_n.\]
Thus, we can construct a sequence of continuous policies which converges to $\gamma$ under the Young topology (at $\psi$).


{\bf Step 3.}
Let $d_{\mathbb{X}}$ denote the metric on $\mathbb{X}$. First assume that $\mathbb{X}$ is compact. For each $m\geq1$, there exists a finite subset $\{z_{m,i}\}_{i=1}^{k_m}$ of $\mathbb{X}$ such that
\begin{align}
\min_{i\in\{1,\ldots,k_m\}} d_{\mathbb{X}}(z,z_{m,i}) < 1/m \text{ for all } z \in \mathbb{X}. \nonumber
\end{align}
Let $\mathbb{X}_m := \{x_{m,1},\ldots,x_{m,k_m}\}$ and define $Q_m$ mapping any $z\in\mathbb{X}$ to the nearest element of $\mathbb{X}_m$, i.e.,
\begin{align}\label{nearestNeighbor}
Q_m(z) := \mathrm{argmin}_{z_{m,i} \in \mathbb{X}_m} d_{\mathbb{X}}(z,z_{m,i}).\nonumber
\end{align}

For each $m$, a partition $\{\mathbb{S}_{m,i}\}_{i=1}^{k_m}$ of the state space $\mathbb{X}$ is induced by $Q_m$ by setting
\begin{align}
\mathbb{S}_{m,i} = \{z \in \mathbb{X}: Q_m(z)=z_{m,i}\}. \nonumber
\end{align}


For a continuous policy $\gamma_C$, the integral
\[\int g(x,u) \gamma_C(du|x),\]
is continuous in $x$: This follows by a generalized convergence theorem since $g(x_n,u) \to g(x,u)$ and $\gamma_C(du|x_n) \to \gamma_C(du|x)$ as as $x_n \to x$; see \cite[Theorem 3.5]{serfozo1982convergence} or \cite[Theorem 3.5]{Lan81}.

Therefore, once we have a continuous policy $\gamma_C$, we can have that for every continuous and bounded $g: \mathbb{X} \times \mathbb{U} \to \mathbb{R}$,
\[\int g(x,u) \gamma_C(du|x) \psi(dx)\]
can be approximated with
\[\sum_{i=1}^{k_m}  \bigg(\int_{x \in \mathbb{S}_{m,i}} \int g(x,u) \gamma_C(du|z_{m,i}) \psi(dx)\bigg) \]
or
\[\sum_{i=1}^{k_m} \bigg( \int g(z_{m,i},u)  \gamma_C(du|z_{m,i})\bigg)  \psi(\mathbb{S}_{m,i})\]
and which in turn can be approximated arbitrarily well, by also quantizing the action space, with
\[\sum_{i=1}^{k_m}  \bigg( \int g(z_{m,i},Q^M(u))\gamma_C(du|z_{m,i}) \bigg) \psi(\mathbb{S}_{m,i}) \]
or
\[\sum_{i=1}^{k_m}  \bigg(\sum_{a_j \in \Lambda_M} g(z_{m,i},a_j)\gamma_C\bigg((Q^M)^{-1}(a_j) \bigg| z_{m,i}\bigg) \bigg) \psi(\mathbb{S}_{m,i}) \]
where $Q^M: \mathbb{U} \to \Lambda_M \subset \mathbb{U}$ is so that $d(u,Q^M(u)) \leq \frac{1}{M}$ and $M$ sufficiently large. 

For the case with non-compact $\mathbb{X}$, by the finiteness of $\psi$, for every $\delta > 0$, there exists a compact set with $\psi(\mathbb{X} \setminus K_{\delta}) < \delta$. 

Thus, one can approximate any stationary policy with quantized policies under the Young topology (at input measure $\psi$) by taking $m$ and $M$ sufficiently large.

{\bf Step 4.}
The proof is then complete by Theorem \ref{ContInvPolicy}: Note that by (\ref{denk2Inv}) as policies $f_n$ converge to $\gamma$, the invariant measure $\uppi_{f_{n}}(\D{x})f_n(\D{u}|x) $ converges in the $w$-$s$ topology, so that
\[\int \uppi_{f_{n}}(\D{x})f_n(\D{u}|x) c(x,u) \to \int \uppi_{\gamma}(\D{x})\gamma(\D{u}|x) c(x,u) \]
for every $c$ which is continuous in $u$ for each $x$. 

{\bf Step 5.}
For (iii) and (iv), first by \cite[Theorem 4.1]{arapostathis2021optimality}, deterministic policies are dense under the Young topology among the space of randomized stationary policies, where the argument critically builds on the fact that the invariant measures are non-atomic (see e.g. \cite[Theorem 3]{milgrom1985distributional} or \cite{balder1997consequences} \cite[Proposition 2.2]{beiglbock2018denseness},
\cite{lacker2018probabilistic}, \cite{borkar1988probabilistic}, or \cite[Chapter 7]{castaing2004young}) and the steps present in the proof of Theorem \ref{ContInvPolicy} (or Step 4 above). The steps above in (i),(ii),(iii) above then applies in this setup for continuous as well as quantized approximations for deterministic stationary policies.
\qed 
\begin{remark}[On approximations and reinforcement learning]
Note that in the above, quantized policies are near-optimal not only for the weakly continuous case under \hyperlink{H1}{(H1)}, but also in the setup under \hyperlink{H2}{(H2)} where continuity of the kernel in the state variable may not hold, which is less explored; quantized policies still lead to convergence of invariant measures under the total variation-weak sense by (\ref{denk2Inv}), and then this leads to approximation of induced integral costs $\int \uppi_{\gamma}(\D{x})\gamma(\D{u}|x) c(x,u)$ with quantized controls (in both state and action). This motivates the application of reinforcement learning methods for problems of this type by restricting control policies to be finitely many. In particular, with having only finitely many policies, {\it win-stay/lose-shift} type algorithms, such as the ones in \cite{yongacoglu2021satisficing} (adapted to average cost control), are applicable: in such algorithms one responds to a Q-learning exploration period, and revises the policy by randomizing over all (finitely many) policies when $\epsilon$-satisfaction of the current policy does not apply and stays at the current policy when such satisfaction applies, and is guaranteed to converge to near optimality in a suitable probabilistic sense (on average, or almost surely if exploration lengths are increasing). 
%
\end{remark}

\begin{remark}[Related approximation results]
One can also show that the set of quantized {\it deterministic} stationary policies $\Gamma_{SD}$ are dense in the space of randomized stationary policies $\Gamma_S$ endowed with the Young topology at input $\psi$ by a further additional argument (given now that by Theorem \ref{denseQuantizedFiniteAC})(ii) quantized randomized policies with finite range are near optimal) via clustering the support sets of each quantization bin to realize a given quantized randomized policy with a finite range. 
We note that Theorem \ref{denseQuantizedFiniteAC} generalizes \cite[Thm. 4.2]{arapostathis2021optimality} which imposed geometric ergodicity for a more restrictive result. We finally note that \cite[Theorem 3.2]{SaldiLinderYukselTAC14} and \cite[Theorem 4.2]{saldi2014near} had established near optimality of quantized policies (though not quantization of the state space), under \hyperlink{H2}{(H2)} (with slightly more restrictive conditions) and \hyperlink{H1}{(H1)}, respectively; and total variation continuity and uniform ergodicity were imposed in \cite[Theorem 4.14]{SaLiYuSpringer} to lead to both state and action quantization (where \cite[Theorem 4.14]{SaLiYuSpringer} also lead to a finite MDP approximation, whereas our result here only shows finite approximation without an approximate MDP construction).

\end{remark}


\section{Conclusion}
We established relations and equivalence conditions between the Young and Borkar topologies. We also showed that under some conditions for a controlled Markov chain, the invariant measure is continuous on the space of stationary control policies defined by a version of Young topology, which makes the space of policies convex and compact.  We finally presented an approximation result involving continuous or quantized (with uniformly quantized range and domain) stationary policies.


\bibliographystyle{plain}

\begin{thebibliography}{10}

\bibitem{arapostathis2010uniform}
A.~Arapostathis and V.~S. Borkar.
\newblock Uniform recurrence properties of controlled diffusions and
  applications to optimal control.
\newblock {\em SIAM Journal on Control and Optimization}, 48(7):4181--4223,
  2010.

\bibitem{survey}
A.~Arapostathis, V.~S. Borkar, E.~Fernandez-Gaucherand, M.~K. Ghosh, and S.~I.
  Marcus.
\newblock Discrete-time controlled {M}arkov processes with average cost
  criterion: A survey.
\newblock {\em SIAM J. Control and Optimization}, 31:282--344, 1993.

\bibitem{arapostathis2012ergodic}
A.~Arapostathis, V.~S. Borkar, and M.~K. Ghosh.
\newblock {\em Ergodic Control of Diffusion Processes}, volume 143.
\newblock Cambridge University Press, 2012.

\bibitem{arapostathis2021optimality}
A.~Arapostathis and S.~Y{\"u}ksel.
\newblock Convex analytic method revisited: Further optimality results and
  performance of deterministic policies in average cost stochastic control.
\newblock {\em Journal of Mathematical Analysis and Applications},
  517(2):126567, 2023.

\bibitem{balder1997consequences}
E.~J. Balder.
\newblock Consequences of denseness of dirac young measures.
\newblock {\em Journal of Mathematical Analysis and Applications},
  207(2):536--540, 1997.

\bibitem{balder2001}
E.~J. Balder.
\newblock On ws-convergence of product measures.
\newblock {\em Mathematics of Operations Research}, 26(3):494--518, 2001.

\bibitem{balder1988generalized}
E.J. Balder.
\newblock Generalized equilibrium results for games with incomplete
  information.
\newblock {\em Mathematics of Operations Research}, 13(2):265--276, 1988.

\bibitem{bauerle2018optimal}
N.~B{\"a}uerle and D.~Lange.
\newblock Optimal control of partially observable piecewise deterministic
  {M}arkov processes.
\newblock {\em SIAM Journal on Control and Optimization}, 56(2):1441--1462,
  2018.

\bibitem{bauerle2010optimal}
N.~B{\"a}uerle and U.~Rieder.
\newblock Optimal control of piecewise deterministic markov processes with
  finite time horizon.
\newblock {\em Modern Trends in Controlled Stochastic Processes: Theory and
  Applications}, pages 123--143, 2010.

\bibitem{bayraktar2022finite}
E.~Bayraktar, N.~B\"auerle, and A.D. Kara.
\newblock Finite approximations and q learning for mean field type multi agent
  control.
\newblock {\em arXiv preprint arXiv:2211.09633}, 2022.

\bibitem{beiglbock2018denseness}
M.~Beiglb{\"o}ck and D.~Lacker.
\newblock Denseness of adapted processes among causal couplings.
\newblock {\em arXiv}, pages arXiv--1805, 2018.

\bibitem{benevs1971existence}
V.~E. Bene{\v{s}}.
\newblock Existence of optimal stochastic control laws.
\newblock {\em SIAM Journal on Control}, 9(3):446--472, 1971.

\bibitem{BeSh78}
D.~P. Bertsekas and S.~E. Shreve.
\newblock {\em Stochastic optimal control: The discrete time case}.
\newblock Academic Press New York, 1978.

\bibitem{Billingsley}
P.~Billingsley.
\newblock {\em Convergence of Probability Measures}.
\newblock Wiley, New York, 1968.

\bibitem{Bil86}
P.~Billingsley.
\newblock {\em Probability and Measure}.
\newblock Wiley, \newblock New York, 2nd edition, 1986.

\bibitem{bismut1973theorie}
J.-M. Bismut.
\newblock {\em Théorie Probabiliste du Contrôle des Diffusions}, volume 181.
\newblock 1973.

\bibitem{bismut1982partially}
J.-M. Bismut.
\newblock Partially observed diffusions and their control.
\newblock {\em SIAM Journal on Control and Optimization}, 20(2):302--309, 1982.

\bibitem{biswas2015mean}
A.~Biswas.
\newblock Mean field games with ergodic cost for discrete time markov
  processes.
\newblock {\em arXiv preprint arXiv:1510.08968}, 2015.

\bibitem{borkar1988probabilistic}
V.~S. Borkar.
\newblock The probabilistic structure of controlled diffusion processes.
\newblock {\em Acta Applicandae Mathematica}, 11(1):19--48, 1988.

\bibitem{borkar1989topology}
V.~S. Borkar.
\newblock A topology for {M}arkov controls.
\newblock {\em Applied Mathematics and Optimization}, 20(1):55--62, 1989.

\bibitem{BorkarRealization}
V.~S. Borkar.
\newblock White-noise representations in stochastic realization theory.
\newblock {\em SIAM J. on Control and Optimization}, 31:1093--1102, 1993.

\bibitem{Borkar2}
V.~S. Borkar.
\newblock Convex analytic methods in {M}arkov decision processes.
\newblock In {\em {H}andbook of Markov Decision Processes, E. A. Feinberg, A.
  Shwartz (Eds.)}, pages 347--375. Kluwer, Boston, MA, 2001.

\bibitem{borkarghosh1990controlled}
V.S. Borkar and M.K. Ghosh.
\newblock Controlled diffusions with constraints.
\newblock {\em Journal of mathematical analysis and applications},
  152(1):88--108, 1990.

\bibitem{castaing2004young}
C.~Castaing, P.~R.~De Fitte, and M.~Valadier.
\newblock {\em Young measures on topological spaces: with applications in
  control theory and probability theory}, volume 571.
\newblock Springer Science \& Business Media, 2004.

\bibitem{CeMe97}
P.~Cembranos and J.~Mendoza.
\newblock {\em Banach Spaces of Vector-Valued Functions}.
\newblock Springer-Verlag, 1997.

\bibitem{dieudonne1951theoreme}
J.~Dieudonn{\'e}.
\newblock Sur le th{\'e}or{\`e}me de lebesgue-nikodym (v).
\newblock {\em Canadian Journal of Mathematics}, 3:129--139, 1951.

\bibitem{Dud02}
R.~M. Dudley.
\newblock {\em Real Analysis and Probability}.
\newblock Cambridge University Press, Cambridge, 2nd edition, 2002.

\bibitem{dugundji}
J.~Dugundji.
\newblock An extension of {T}ietze's theorem.
\newblock {\em Pacific Journal of Mathematics}, 1(3):353--367, 1951.

\bibitem{duncan1971solutions}
T.~Duncan and P.~Varaiya.
\newblock On the solutions of a stochastic control system.
\newblock {\em SIAM Journal on Control}, 9(3):354--371, 1971.

\bibitem{ethier2009markov}
S.~N. Ethier and T.~G. Kurtz.
\newblock {\em Markov Processes: Characterization and Convergence}, volume 282.
\newblock John Wiley \& Sons, 2009.

\bibitem{fattorini1994existence}
H.O. Fattorini.
\newblock Existence theory and the maximum principle for relaxed
  infinite-dimensional optimal control problems.
\newblock {\em SIAM Journal on Control and Optimization}, 32(2):311--331, 1994.

\bibitem{fleming1976generalized}
W.~H. Fleming.
\newblock Generalized solutions in optimal stochastic control.
\newblock In {\em Proc. U.R.I. Conference on Control, University of Rhode
  Island}, 1982.

\bibitem{florescu2012young}
L.C. Florescu and C.~Godet-Thobie.
\newblock {\em Young measures and compactness in measure spaces}.
\newblock Walter de Gruyter, 2012.

\bibitem{Fol99}
G.~B. Folland.
\newblock {\em Real Analysis: Modern Techniques and Their Applications}.
\newblock John Wiley and Sons, 1999.

\bibitem{hernandez1993existence}
O.~Hern{\'a}ndez-Lerma.
\newblock Existence of average optimal policies in markov control processes
  with strictly unbounded costs.
\newblock {\em Kybernetika}, 29(1):1--17, 1993.

\bibitem{HernandezLermaMCP}
O.~Hern{\'a}ndez-Lerma and J.~B. Lasserre.
\newblock {\em Discrete-Time {M}arkov Control Processes: Basic Optimality
  Criteria}.
\newblock Springer, 1996.

\bibitem{HernandezLasserreErgodic}
O.~Hern{\'a}ndez-Lerma and J.~B. Lasserre.
\newblock {\em {M}arkov Chains and Invariant Probabilities}.
\newblock Birkh\"auser, Basel, 2003.

\bibitem{KYSICONPrior}
A.D Kara and S.~Y\"uksel.
\newblock Robustness to incorrect priors in partially observed stochastic
  control.
\newblock {\em SIAM Journal on Control and Optimization}, 57(3):1929--1964,
  2019.

\bibitem{kartashov1996strong}
N.V. Kartashov.
\newblock {\em Strong stable Markov chains}.
\newblock De Gruyter, 1996.

\bibitem{kurano1989existence}
M.~Kurano.
\newblock The existence of a minimum pair of state and policy for markov
  decision processes under the hypothesis of doeblin.
\newblock {\em SIAM journal on control and optimization}, 27(2):296--307, 1989.

\bibitem{kushner2014partial}
H.~J. Kushner.
\newblock A partial history of the early development of continuous-time
  nonlinear stochastic systems theory.
\newblock {\em Automatica}, 50(2):303--334, 2014.

\bibitem{kushner2001numerical}
H.~J. Kushner and P.~G. Dupuis.
\newblock {\em Numerical Methods for Stochastic Control Problems in Continuous
  Time}, volume~24.
\newblock Springer Science \& Business Media, 2001.

\bibitem{lacker2018probabilistic}
D.~Lacker.
\newblock Probabilistic compactification methods for stochastic optimal control
  and mean field games.
\newblock 2018.

\bibitem{Lan81}
H.J. Langen.
\newblock Convergence of dynamic programming models.
\newblock {\em Mathematics of Operations Research}, 6(4):493--512, Nov. 1981.

\bibitem{mascolo1989relaxation}
E.~Mascolo and L.~Migliaccio.
\newblock Relaxation methods in control theory.
\newblock {\em Applied Mathematics and Optimization}, 20(1):97--103, 1989.

\bibitem{mcshane1967relaxed}
E.~J. McShane.
\newblock Relaxed controls and variational problems.
\newblock {\em SIAM Journal on Control}, 5(3):438--485, 1967.

\bibitem{mertens2015repeated}
J.-F. Mertens, S.~Sorin, and S.~Zamir.
\newblock {\em Repeated games}, volume~55.
\newblock Cambridge University Press, 2015.

\bibitem{milgrom1985distributional}
P.~R. Milgrom and R.~J. Weber.
\newblock Distributional strategies for games with incomplete information.
\newblock {\em Mathematics of operations research}, 10(4):619--632, 1985.

\bibitem{Neveu}
J.~Neveu.
\newblock {\em Mathematical foundations of the calculus of probability}.
\newblock Holden-Day, Inc., San Francisco, Calif.-London-Amsterdam, 1965.

\bibitem{Choquet}
R.P. Phelps.
\newblock {\em Lectures on Choquet's theorem}.
\newblock Van Nostrand, New York:, 1966.

\bibitem{pradhan2022near}
S.~Pradhan and S.~Y{\"u}ksel.
\newblock Continuity of cost in {B}orkar control topology and implications on
  discrete space and time approximations for controlled diffusions under
  several criteria.
\newblock {\em arXiv preprint arXiv:2209.14982}, 2022.

\bibitem{rudin1991functional}
W.~Rudin.
\newblock {\em Functional Analysis}.
\newblock McGraw-Hill, 1991.

\bibitem{SaldiArXiv2017}
N.~Saldi.
\newblock A topology for team policies and existence of optimal team policies
  in stochastic team theory.
\newblock {\em IEEE Transactions on Automatic Control}, 65(1):310--317, 2020.

\bibitem{SaldiLinderYukselTAC14}
N.~Saldi, T.~Linder, and S.~Y\"uksel.
\newblock Asymptotic optimality and rates of convergence of quantized
  stationary policies in stochastic control.
\newblock {\em IEEE Trans. Automatic Control}, 60:553 --558, 2015.

\bibitem{SaLiYuSpringer}
N.~Saldi, T.~Linder, and S.~Y\"uksel.
\newblock {\em Finite Approximations in Discrete-Time Stochastic Control:
  Quantized Models and Asymptotic Optimality}.
\newblock Springer, Cham, 2018.

\bibitem{saldiyukselGeoInfoStructure}
N.~Saldi and S.~Y\"uksel.
\newblock Geometry of information structures, strategic measures and associated
  control topologies.
\newblock {\em Probability Surveys}, 19:450--532, 2022.

\bibitem{saldi2014near}
N.~Saldi, S.~Y\"uksel, and T.~Linder.
\newblock Near optimality of quantized policies in stochastic control under
  weak continuity conditions.
\newblock {\em Journal of Mathematical Analysis and Applications},
  435(1):321--337, 2016.

\bibitem{schal1975dynamic}
M.~Sch{\"a}l.
\newblock On dynamic programming: compactness of the space of policies.
\newblock {\em Stochastic Processes and their Applications}, 3(4):345--364,
  1975.

\bibitem{serfozo1982convergence}
R.~Serfozo.
\newblock Convergence of {L}ebesgue integrals with varying measures.
\newblock {\em Sankhy{\=a}: The Indian Journal of Statistics, Series A}, pages
  380--402, 1982.

\bibitem{warga2014optimal}
J.~Warga.
\newblock {\em Optimal Control of Differential and Functional Equations}.
\newblock Academic press, 2014.

\bibitem{yongacoglu2021satisficing}
B.~Yongacoglu, G.~Arslan, and S.~Y\"uksel.
\newblock Satisficing paths and independent multi-agent reinforcement learning
  in stochastic games.
\newblock {\em SIAM Journal on Mathematics of Data Science (arXiv:2110.04638)},
  2023.

\bibitem{young1937generalized}
L.C. Young.
\newblock Generalized curves and the existence of an attained absolute minimum
  in the calculus of variations.
\newblock {\em Comptes Rendus de la Societe des Sci. et des Lettres de
  Varsovie}, 30:212--234, 1937.

\bibitem{YukselWitsenStandardArXiv}
S.~Y\"uksel.
\newblock A universal dynamic program and refined existence results for
  decentralized stochastic control.
\newblock {\em SIAM Journal on Control and Optimization}, 58(5):2711--2739,
  2020.

\bibitem{YukselOptimizationofChannels}
S.~Y\"uksel and T.~Linder.
\newblock Optimization and convergence of observation channels in stochastic
  control.
\newblock {\em SIAM J. on Control and Optimization}, 50:864--887, 2012.

\end{thebibliography}

\end{document}